\documentclass[10pt,twoside]{article}
\usepackage{amssymb}

{\catcode`\@=11
\gdef\ps@myheadings{\let\@mkboth\@gobbletwo
 \def\@oddhead
          {\vbox{\noindent
                {\scriptsize\bf Key Renewal Theorem
                                         \hfill\thepage}\vskip 4pt}}%
 \def\@oddfoot{}
 \def\@evenhead
          {\vbox{\noindent{\scriptsize\bf\thepage\hfill
          Korshunov}\vskip 4pt}}%
 \def\@evenfoot{}\def\chaptermark##1{}\def\sectionmark##1{}%
 \def\subsectionmark##1{}}

\gdef\@begintheorem#1#2{\trivlist \item[\hskip \labelsep{\indent\bf #1\ #2.}]\it}
\gdef\@opargbegintheorem#1#2#3{\trivlist
      \item[\hskip \labelsep{\bf #1\ #2\ (#3)}]\it}
\gdef\@endtheorem{\endtrivlist}
 }

\newcommand\proof{P\,r\,o\,o\,f}
\newtheorem{Theorem}{Theorem}
\newtheorem{Lemma}{Lemma}
\newtheorem{Corollary}{Corollary}

\begin{document}
\thispagestyle{empty}
\section*{\large\bf
           The Key Renewal Theorem
           for a Transient Markov Chain}
\section*{\normalsize\bfseries
Dmitry Korshunov\footnote{Sobolev Institute of Mathematics,
                      Novosibirsk 630090, Russia.\\
                      E-mail: korshunov@math.nsc.ru}}

\begin{abstract}
\noindent
{\scriptsize We consider a time-homogeneous
Markov chain $X_n$, $n\ge0$, valued in ${\bf R}$.
Suppose that this chain is transient, that is,
$X_n$ generates a $\sigma$-finite renewal measure.
We prove the key renewal theorem under condition
that this chain has asymptotically homogeneous
at infinity jumps and asymptotically positive drift.

\noindent
{\bf KEY WORDS}:
transient Markov chain,
renewal kernel,
renewal measure,
key renewal theorem,
Green function.}
\end{abstract}

\noindent
Let $\xi_1$, $\xi_2$, \ldots\ be independent
identically distributed random variables
with a common distribution $F$ on ${\bf R}$.
Put $S_n=\xi_1+\cdots+\xi_n$, $S_0=0$,
and consider the renewal measure generated by sums:
\begin{eqnarray*}
U(B) &\equiv& {\bf E}\sum_{n=0}^\infty {\bf I}\{S_n\in B\}
=\sum_{n=0}^\infty F^{*n}(B).
\end{eqnarray*}
If $F$ is non-lattice then
the celebrated key renewal theorem states that,
for every fixed $h>0$,
\begin{eqnarray}\label{renewal.sums}
U(x,x+h] &\to& \frac{h}{{\bf E}\xi_1}\ \mbox{ as }x\to\infty,
\end{eqnarray}
provided ${\bf E}\xi_1$ is finite and positive
(see, for example, Feller and Orey [\ref{FO}],
Feller [\ref{F}, Ch. XI],
Woodroofe [\ref{Woodroofe}, Appendix]);
$F$ is called lattice if it is concentrated
on some lattice $\{ka,k\in{\bf Z}\}$ with $a>0$.
If $F$ is lattice the same is true when $h$
is a multiple of the span $a$.

It is proved in Wang and Woodroofe [\ref{WW}] and
in Borovkov and Foss [\ref{BF}, Theorem 2.6]
that (\ref{renewal.sums}) holds uniformly over
certain classes of distributions $F$.
Some extensions of the key renewal theorem
to the nonidentically distributed case
are considered by Williamson [\ref{W1965}],
Maejima [\ref{M1975}].
Another extension is aimed to include random walks
perturbed by both a slowly changing sequence
and a stationary one, $Z_n=S_n+\eta_n+\zeta_n$ say;
see, for example, Lai and Siegmund [\ref{LS1}, \ref{LS2}],
Woodroofe [\ref{Woodroofe}], Horv\'ath [\ref{Horvath}],
Zhang [\ref{Zhang}], Kim and Woodroofe [\ref{KW}].
All these extensions deal with perturbations depending
on time rather than on state space; in most cases
summands are independent and have finite variance.

Many papers (see, for example, Kesten [\ref{Kesten}],
Athreya, McDonald and Ney [\ref{ADN}],
Nummelin [\ref{Nummelin}], Alsmeyer [\ref{Alsmeyer}],
Kl\"uppelberg and Pergamenchtchikov [\ref{KP}],
Fuh [\ref{Fuh}], and also some of their references)
are devoted to a Markov modulated random walks.
Usually a (Harris-) recurrent Markov chain $X_n$
is considered with an invariant measure, $\pi$ say.
Conditioned on a realisation $\{x_n,n\ge0\}$,
one is given a sequence of independent random variables
$\{\xi_n\}$, such that the distribution of $\xi_n$
depends only on $x_n$. Put $T_n=\xi_0+\ldots+\xi_n$
and assume $T_n/n\to\alpha$ with probability $1$.
Then the typical result states the convergence
\begin{eqnarray*}
{\bf E}\sum_{n=0}^\infty g(X_n,x-T_n)
&\to& \frac{1}{\alpha}\int\pi(dy)\int_{-\infty}^\infty g(y,s)ds
\quad\mbox{ as }x\to\infty,
\end{eqnarray*}
for bounded continuous function $g$ satisfying some
conditions.
The corresponding proofs use probabilistic arguments,
notably the construction of regeneration epochs
for $\{X_n,n\ge0\}$ (for example, visit times
to some atom). This approach eventually reduces
the problem to Blackwell's renewal theorem for sums
of independent identically distributed random variables.

To the best of our knowledge, the only result related
to a random walk perturbed in state space is due to
Heyde [\ref{Heyde}]. To be more precise,
Heyde discussed the key renewal theorem for maxima
$M_n\equiv\max\limits_{0\le k\le n}S_k$ of partial sums:
provided ${\bf E}\xi_1$ is finite and positive
\begin{eqnarray}\label{renewal.max}
{\bf E}\sum_{n=0}^\infty {\bf I}\{M_n\in(x,x+h]\}
&\to& \frac{h}{{\bf E}\xi_1}\ \mbox{ as }\ x\to\infty.
\end{eqnarray}
It is well known (see, for example, [\ref{F}, Ch. VI, Sec. 9])
that $M_n$ has the same distribution as the reflected
random walk on ${\bf R}^+$ defined by the recursion
\begin{eqnarray*}
W_{n+1} &=& (W_n+\xi_{n+1})^+, \quad W_0=0.
\end{eqnarray*}
So, the renewal function generated
by the chain $W_n$ has the same asymptotic behaviour
as described in (\ref{renewal.max}).

Both the random walk $S_n$ and the reflected random walk
on the positive half-line $W_n$ are particular
examples of Markov chains on ${\bf R}$.
In the present paper we extend
the key renewal theorem from these
very important cases onto asymptotically
space-homogeneous Markov chains on ${\bf R}$
with an asymptotically positive drift.
Introduce some relevant definitions.

Let $P(x,B)$, $x\in{\bf R}$, $B\in{\mathcal B}({\bf R})$,
be a transition probability kernel on ${\bf R}$;
hereinafter ${\mathcal B}({\bf R})$
is the Borel $\sigma$-algebra on ${\bf R}$.
Consider a time-homogeneous Markov chain
$X=\{X_n,\ n=0,1,2,\ldots\}$ on ${\bf R}$
with transition probabilities $P(\cdot,\cdot)$, that is,
\begin{eqnarray*}
{\bf P}\{X_{n+1}\in B\,|\,X_n=x\} &=& P(x,B).
\end{eqnarray*}
Let $\xi(x)$ be the random variable distributed
as the jump of the chain at state $x$:
\begin{eqnarray*}
{\bf P}\{x+\xi(x)\in B\} &=& P(x,B),
\quad B\in{\mathcal B}({\bf R}).
\end{eqnarray*}
Let $\mu_n$ denote the distribution of $X_n$;
then equalities $\mu_n=\mu_{n-1}P$
and $\mu_n=\mu_0 P^n$ hold.
Formally, define the renewal (or potential)
kernel $Q$ by the equality
\begin{eqnarray*}
Q(\cdot,\cdot) &=& \sum_{n=0}^\infty P^n(\cdot,\cdot).
\end{eqnarray*}
We assume that the Markov chain $X$ is
{\it transient} (see Meyn and Tweedie [\ref{Tw}, Ch. 8]),
that is, there exists a countable cover of
${\bf R}$ with uniformly transient sets $\{B_k\}$.
In its turn a set $B\in{\mathcal B}({\bf R})$
is called {\it uniformly transient} if
\begin{eqnarray}\label{uniform.reccurence}
\sup_{y\in B}Q(y,B) &<& \infty.
\end{eqnarray}
By the Markov property, this is equivalent to
\begin{eqnarray}\label{uniform.reccurence.whole}
\sup_{y\in{\bf R}}Q(y,B) &<& \infty.
\end{eqnarray}
Indeed, considering the first hitting time of $B$,
we conclude the following inequality,
for each state $x\in{\bf R}$,
\begin{eqnarray}\label{recc.whole.B}
Q(x,B) &\le& \sup_{y\in B}Q(y,B),
\end{eqnarray}
which implies (\ref{uniform.reccurence.whole}).

In the present paper we assume that $X_n$
is transient with respect to the collection of sets
$B_k=(k,k+1]$, $k\in{\bf Z}$; that is,
for any $k\in{\bf Z}$,
\begin{eqnarray}\label{uniform.reccurence.our}
\sup_{y\in{\bf R}}Q(y,(k,k+1]) &<& \infty.
\end{eqnarray}
Then $Q(x,B)<\infty$ for all $x$ and bounded set $B$.
Hence the renewal measure generated by the chain $X$
\begin{eqnarray*}
U(B) &\equiv& \sum_{n=0}^\infty {\bf P}\{X_n\in B\}
= \sum_{n=0}^\infty \mu_n(B) = (\mu_0 Q)(B)
\end{eqnarray*}
is well defined for every initial distribution
$\mu_0$ and bounded set $B$;
$U$ is $\sigma$-finite with respect to the collection
of sets $(k,k+1]$, $k\in{\bf Z}$.

The main goal of our analysis is the local
asymptotic behaviour of this renewal measure.
Without further restrictions on the chain $X_n$,
the asymptotics of $U(x,x+h]$ as $x\to\infty$
can be very special.
We consider a transient Markov chain $X_n$ as
{\it a perturbation in space} of the random walk $S_n$
with positive drift.
To get similar renewal behaviour
for $X_n$ as for $S_n$, it is natural to assume that,
being far away from the origin,
$X_n$ behaves almost like $S_n$.

Thus we restrict our attention to the
{\it asymptotically space-ho\-mo\-ge\-neous
Markov chain} $X$, that is, we assume that the
distribution of the jump $\xi(x)$ has a weak limit $F$
as $x\to\infty$. Let $\xi$ be a random variable
with distribution $F$.

The notion of asymptotically space-homogeneous
Markov chain is a natural generalisation of both
(i) the random walk $S_n$;
in this case $\xi(x)=_{\rm d}\xi$ for all $x$;
(ii) the reflected random walk $W_n$ on the positive half-line;
in this case $\xi(x)=_{\rm d}(x+\xi)^+-x$.
An asymptotically space-homogeneous Markov chains appear
in different areas; in particular, we are motivated by
theory of queues when the service rate depends
on the current waiting time;
and by sequential analysis related to an optimal solutions
in a change-point problem (see Borovkov [\ref{Borovkov}]).
Some limit theorems for them were obtained by
Korshunov [\ref{Korshunov}].

So, let the Markov chain $X_n$ be
asymptotically space-homogeneous.

\begin{Theorem}\label{renewal.theorem}
Let $\xi(x)\Rightarrow\xi$ as $x\to\infty$
and ${\bf E}\xi>0$.
Let the family of random variables
$\{|\xi(x)|,\ x\in{\bf R}\}$
admit an integrable majorant $\eta$,
that is, ${\bf E}\eta<\infty$ and
\begin{eqnarray}\label{majoriz}
|\xi(x)| &\le_{\rm st}& \eta
\quad\mbox{for all }x\in{\bf R}.
\end{eqnarray}
Assume that
\begin{eqnarray}\label{finite.bound.U}
\sup_{k\in{\bf Z}} U(k,k+1]
&<& \infty.
\end{eqnarray}
Assume also that there exists a limit
\begin{eqnarray}\label{lim.Xge0}
p_0 &=& \lim_{n\to\infty}{\bf P}\{X_n>0\}.
\end{eqnarray}

If the limit distribution $F$ is non-lattice,
then $U(x,x+h]\to h/{\bf E}\xi$
as $x\to\infty$, for every fixed $h>0$.

If the chain $X_n$ is integer valued and
${\bf Z}$ is the lattice with minimal span
for distribution $F$, then
$U\{n\}\to 1/{\bf E}\xi$ as $n\to\infty$.
\end{Theorem}

Condition (\ref{majoriz}) and
the dominated convergence theorem imply
$|\xi|\le_{\rm st}\eta$, ${\bf E}|\xi|<\infty$
and ${\bf E}\xi(x)\to{\bf E}\xi$ as $x\to\infty$;
in particular, the chain $X_n$
has an asymptotically space-homogeneous drift.

Only conditions (\ref{finite.bound.U}) and
(\ref{lim.Xge0}) of Theorem \ref{renewal.theorem}
are not formulated in local terms, i.e.,
in terms of one-step transition probabilities.
Below, in Theorem \ref{suff.for.tran},
we give some simple conditions sufficient
for (\ref{finite.bound.U}).
Note that the value of $p_0$ in condition
(\ref{lim.Xge0}) may be very sensitive
with respect to the local probabilities.
It can be illustrated by the following example.
Let $X_n$ be a chain valued on ${\bf Z}$ with
the following transition probabilities:
\begin{eqnarray*}
p_{i,i+1} = 3/4, && p_{i,i-1}=1/4
\quad\mbox{ for } i\ge 1,\\
p_{i,i+1} = 1/4, && p_{i,i-1}=3/4
\quad\mbox{ for } i\le -1,\\
p_{0,1} = p, && p_{0,-1}=1-p,
\end{eqnarray*}
where $p\in[0,1]$. Given $X_0=0$, then
$p_0=p_0(p)=p$ is increasing from $0$ to $1$
simultaneously with $p$.

Since the chain is transient,
by (\ref{uniform.reccurence.whole})
the convergence
$\mu_n(K)\equiv{\bf P}\{X_n\in K\}\to 0$
holds as $n\to\infty$ for any compact $K$.
Hence, condition (\ref{lim.Xge0}) is equivalent
to the convergence, for every fixed $x_0$,
\begin{eqnarray}\label{lim.Xge1}
{\bf P}\{X_n>x_0\} &\to& p_0
\quad \mbox{ as }n\to\infty.
\end{eqnarray}

\proof\ of Theorem \ref{renewal.theorem}
follows some ideas of the operator approach
proposed by Feller [\ref{F}, Ch. XI].
First of all, condition (\ref{finite.bound.U})
allows us to apply Helly's Selection Theorem
to the family of measures $\{U(k+\cdot),k\in{\bf Z}^+\}$
(see, for example, Theorem 2 in [\ref{F},
Ch. VIII, Sec. 6]).
Hence, there exists a sequence of points
$t_n\to\infty$ such that the sequence of measures
$U_n(\cdot)\equiv U(t_n+\cdot)$ converges weakly
to some measure $\lambda$ as $n\to\infty$.
The following two lemmas describe properties
of $\lambda$.

\begin{Lemma}\label{l.1}
A weak limit $\lambda$ of the sequence of measures
$U(t_n+\cdot)$ satisfies the identity $\lambda=\lambda*F$.
\end{Lemma}

\proof.
The measure $\lambda$ is non-negative and
$\sigma$-finite with necessity.
Fix any smooth function $f(x)$ with a bounded support;
let $A>0$ be such that $f(x)=0$ for $x\notin[-A,A]$.
The weak convergence of measures means
the convergence of integrals
\begin{eqnarray}\label{conv.f.1}
\int_{-\infty}^\infty f(x)U(t_n+dx)
\equiv \int_{-A}^A f(x)U(t_n+dx)
&\to& \int_{-A}^A f(x)\lambda(dx)
\end{eqnarray}
as $n\to\infty$.
On the other hand, due to the equality $U=\mu_0+UP$
we have the following representation for the
left side of (\ref{conv.f.1}):
\begin{eqnarray}\label{conv.f.2}
\int_{-A}^A f(x)\mu_0(t_n+dx)
+\int_{-A}^A f(x)
\int_{-\infty}^\infty P(t_n+y,t_n+dx)U(t_n+dy).
\end{eqnarray}
Since $f$ is bounded and $\mu_0$ is finite,
\begin{eqnarray}\label{conv.f.3}
\int_{-A}^A f(x)\mu_0(t_n+dx)
&\le& ||f||_C \mu_0[t_n-A,t_n+A] \to 0
\end{eqnarray}
as $n\to\infty$. The second term in (\ref{conv.f.2})
is equal to
\begin{eqnarray}\label{conv.f.4}
\int_{-\infty}^\infty U(t_n+dy)
\int_{-A}^A f(x)P(t_n+y,t_n+dx).
\end{eqnarray}
The weak convergence $P(t,t+\cdot)\Rightarrow F(\cdot)$
as $t\to\infty$ implies the convergence of the inner
integral in (\ref{conv.f.4}):
\begin{eqnarray*}
\int_{-A}^A f(x)P(t_n+y,t_n+dx)
&\to& \int_{-A}^A f(x)F(dx-y);
\end{eqnarray*}
here the rate of convergence can be estimated
in the following way:
\begin{eqnarray*}
\Delta(n,y) &\equiv& \Biggl|\int_{-A}^A
f(x) (P(t_n+y,t_n+dx)-F(dx-y))\Biggr|\\
&=& \Biggl|\int_{-A}^A
f'(x)({\bf P}\{\xi(t_n+y)\le x-y\}-F(x-y))dx\Biggr|\\
&\le& ||f'||_C \int_{-A-y}^{A-y}
|{\bf P}\{\xi(t_n+y)\le x\}-F(x)|dx.
\end{eqnarray*}
Thus, the asymptotic homogenuity of the chain
yields for every fixed $C>0$ the uniform convergence
\begin{eqnarray}\label{Delta.1}
\sup_{y\in[-C,C]}\Delta(n,y) &\to& 0
\quad\mbox{as }n\to\infty.
\end{eqnarray}
In addition, by majorisation condition
(\ref{majoriz}), for all $x$
regardless positive or negative,
\begin{eqnarray*}
|{\bf P}\{\xi(t_n+y)\le x\}-F(x)|
&\le& 2{\bf P}\{\eta>|x|\}.
\end{eqnarray*}
Hence, for all $y$,
\begin{eqnarray}\label{Delta.2}
\Delta(n,y) &\le& 2||f'||_C \int_{-A-y}^{A-y}
{\bf P}\{\eta>|x|\}dx\nonumber\\
&\le& 4A||f'||_C {\bf P}\{\eta>|y|-A\}.
\end{eqnarray}
We have the estimate
\begin{eqnarray*}
\Delta_n &\equiv&
\Biggl|\int_{-\infty}^\infty U(t_n+dy)
\Biggl(\int_{-\infty}^\infty f(x)P(t_n{+}y,t_n{+}dx)
-\int_{-\infty}^\infty f(x)F(dx{-}y)\Biggr)\Biggr|\\
&\le& \int_{-\infty}^\infty \Delta(y,n) U(t_n+dy).
\end{eqnarray*}
For any fixed $C>0$, uniform convergence
(\ref{Delta.1}) implies
\begin{eqnarray*}
\int_{-C}^C \Delta(y,n) U(t_n+dy)
&\le& \sup_{y\in[-C,C]} \Delta(y,n)
\cdot\sup_n U[t_n-C,t_n+C]\\
&\to& 0 \quad\mbox{as }n\to\infty.
\end{eqnarray*}
The remaining part of the integral can be
estimated by (\ref{Delta.2}):
\begin{eqnarray*}
\lefteqn{\limsup_{n\to\infty}\int_{|y|\ge C} \Delta(y,n) U(t_n+dy)}\\
&\le& 4A||f'||_C \limsup_{n\to\infty}
\int_{|y|\ge C} {\bf P}\{\eta>|y|-A\} U(t_n+dy).
\end{eqnarray*}
Since $\eta$ has a finite mean, property
(\ref{finite.bound.U}) of the renewal measure $U$
allows us to choose a sufficiently large $C$
in order to make the `$\limsup$' as small as we please.
Therefore, $\Delta_n \to 0$ as $n\to\infty$.
Hence, (\ref{conv.f.4}) has the same limit
as the sequence of integrals
\begin{eqnarray*}
\int_{-\infty}^\infty U(t_n+dy)
\int_{-A}^A f(x)F(dx-y).
\end{eqnarray*}
Now the weak convergence to $\lambda$
implies that (\ref{conv.f.4}) has the limit
\begin{eqnarray}\label{conv.f.5}
\int_{-\infty}^\infty \lambda(dy)
\int_{-\infty}^\infty f(x)F(dx-y)
&=& \int_{-\infty}^\infty f(x)
\int_{-\infty}^\infty F(dx-y) \lambda(dy)\nonumber\\
&=& \int_{-\infty}^\infty f(x) (F*\lambda)(dx).
\end{eqnarray}
By (\ref{conv.f.1})--(\ref{conv.f.3})
and (\ref{conv.f.5}), we conclude the identity
\begin{eqnarray*}
\int_{-\infty}^\infty f(x)\lambda(dx)
&=& \int_{-\infty}^\infty f(x) (F*\lambda)(dx).
\end{eqnarray*}
Since this identity holds for every smooth function
$f$ with a bounded support, the measures
$\lambda$ and $F*\lambda$ coincide.
The proof is complete.

Further we use the following statement
which was proved in [\ref{CD}]
(see also [\ref{Raugi}] or [\ref{Revuz}, Sec. 5.1]):

\begin{Lemma}\label{l.2}
Let $F$ be a distribution not concentrated at $0$.
Let $\lambda$ be a nonnegative measure satisfying
the equality $\lambda=F*\lambda$ and the property
$\sup\limits_{n\in{\bf Z}}\lambda[n,n+1]<\infty$.

If $F$ is non-lattice,
then $\lambda$ is proportional to Lebesgue measure.

If $F$ is lattice with minimal span $1$ and
$\lambda({\bf Z})=1$, then $\lambda$ is
proportional to the counting measure.
\end{Lemma}


The concluding part of the proof of
Theorem \ref{renewal.theorem}
will be carried out for the non-lattice case.
Choose any sequence of points $t_n\to\infty$
such that the measure $U(t_n+\cdot)$ converges
weakly to some measure $\lambda$ as $n\to\infty$.
It follows from Lemmas \ref{l.1} and \ref{l.2}
that then $\lambda(dx)=\alpha\cdot dx$
with some $\alpha$, i.e.,
\begin{eqnarray*}
U(t_n+dx) &\Rightarrow& \alpha\cdot dx
\ \mbox{ as }n\to\infty.
\end{eqnarray*}
Now it suffices to prove that $\alpha=p_0/{\bf E}\xi$.

Fix some $k\in{\bf N}$. Put
$U_k\equiv UP^k=\sum_{j=k}^\infty\mu_j$. Then
\begin{eqnarray}\label{Uk.tn.to.lambda}
U_k(t_n+dx) &\Rightarrow& \alpha\cdot dx
\ \mbox{ as }n\to\infty.
\end{eqnarray}
Consider the measure $U_k-U_{k+1}=U_k(I-P)$;
by the definition of the renewal measure
it is equal to $\mu_k$, that is,
for any {\it bounded} Borel set $B$,
$U_k(B)-U_{k+1}(B)=\mu_k(B)$
(the equality may fail for unbounded sets,
say, for $(-\infty,x]$).
In particular,
\begin{eqnarray}\label{I.-.P.to.1}
(U_k-U_{k+1})(0,x] &=& \mu_k(0,x]
\to \mu_k(0,\infty)
\ \mbox{ as }x\to\infty.
\end{eqnarray}
On the other hand,
\begin{eqnarray}\label{I.-.P.+.-}
\lefteqn{(U_k-U_{k+1})(0,x]
= \int_{-\infty}^\infty (I-P)(y,(0,x])U_k(dy)\nonumber}\\
&=& - \int_{-\infty}^0 P(y,(0,x])U_k(dy)
+\int_0^x P(y,(-\infty,0])U_k(dy)\nonumber\\
&& +\int_0^x P(y,(x,\infty))U_k(dy)
-\int_x^\infty P(y,(0,x])U_k(dy).
\end{eqnarray}
The asymptotic homogeneity of the chain and
weak convergence (\ref{Uk.tn.to.lambda})
imply the following convergences
of the integrals, for any fixed $A>0$:
\begin{eqnarray}\label{I.-.P.+.-.1}
\int_{t_n-A}^{t_n} P(y,(t_n,\infty))U_k(dy)
&\to& \alpha \int_0^A {\bf P}\{\xi>z\}dz
\end{eqnarray}
as $n\to\infty$, and
\begin{eqnarray}\label{I.-.P.+.-.2}
\int_{t_n}^{t_n+A} P(y,(0,t_n])U_k(dy)
&\to& \alpha \int_0^A {\bf P}\{\xi\le -z\}dz.
\end{eqnarray}
Majorisation condition (\ref{majoriz})
allows us to estimate the tails of the integrals:
\begin{eqnarray}\label{I.-.P.+.-.3}
\int_0^{t_n-A} P(y,(t_n,\infty))U_k(dy)
&\le& -\int_A^\infty {\bf P}\{\eta>z\}U(t_n-dz)
\end{eqnarray}
and
\begin{eqnarray}\label{I.-.P.+.-.4}
\int_{t_n+A}^\infty P(y,(0,t_n])U_k(dy)
&\le& \int_A^\infty {\bf P}\{\eta\ge z\}U(t_n+dz).
\end{eqnarray}
Since the majorant $\eta$ is integrable,
condition (\ref{finite.bound.U})
guarantees that the right sides of inequalities
(\ref{I.-.P.+.-.3}) and (\ref{I.-.P.+.-.4})
can be made as small as we please by the choice
of sufficiently large $A$.
By these reasons we conclude from
(\ref{I.-.P.+.-})--(\ref{I.-.P.+.-.2}) that,
as $n\to\infty$,
\begin{eqnarray*}
\lefteqn{(U_k-U_{k+1})(0,t_n]}\\
&\to& -\int_{-\infty}^0 P(y,(0,\infty))U_k(dy)
+\int_0^\infty P(y,(-\infty,0])U_k(dy)\\
&&\hspace{40mm} +\alpha\int_0^\infty {\bf P}\{\xi>z\}dz
-\alpha\int_0^\infty {\bf P}\{\xi\le -z\}dz.
\end{eqnarray*}
Together with (\ref{I.-.P.to.1})
it implies the following equality,
for any fixed $k$:
\begin{eqnarray}\label{mu.k.int0.alpha}
\mu_k(0,\infty)
&=& -\int_{-\infty}^0 P(y,(0,\infty))U_k(dy)
+\int_0^\infty P(y,(-\infty,0])U_k(dy)
+\alpha{\bf E}\xi.\nonumber\\[-2mm]
\end{eqnarray}
Now let $k\to\infty$, then both integrals go to zero.
For example, the first integral can be estimated
in the following way, for every $A>0$:
\begin{eqnarray*}
\int_{-\infty}^0 P(y,(0,\infty))U_k(dy)
&\le& \int_{-\infty}^{-A} {\bf P}\{\eta>-y\}U(dy)
+U_k(-A,0].
\end{eqnarray*}
Here, for any fixed $A$,
$U_k(-A,0]\to 0$ as $k\to\infty$
by (\ref{uniform.reccurence.whole}).
Therefore, it follows from (\ref{mu.k.int0.alpha})
and (\ref{lim.Xge0}) that
$p_0=\alpha{\bf E}\xi$. The proof
of Theorem \ref{renewal.theorem} is complete.

In the next theorem we provide some simple
conditions sufficient for condition (\ref{finite.bound.U}),
that is, for local compactness of the renewal measure.
Denote $a\wedge b=\min\{a,b\}$.

\begin{Theorem}\label{suff.for.tran}
Suppose that there exists $A>0$ such that
\begin{eqnarray}\label{sft.1}
\varepsilon &\equiv&
\inf_{x\in{\bf R}} {\bf E}(\xi(x)\wedge A) > 0.
\end{eqnarray}
In addition, let
\begin{eqnarray}\label{sft.2}
\delta &\equiv& \inf_{x\in{\bf R}}
{\bf P}\{X_n>x\mbox{ for all } n\ge 1|X_0=x\} > 0.
\end{eqnarray}
Then $U(x,x+h]\le (A+h)/\varepsilon\delta$
for all $x\in{\bf R}$ and $h>0$; in particular,
{\rm(\ref{finite.bound.U})} holds.
\end{Theorem}

\proof. Inequality (\ref{recc.whole.B}) implies
\begin{eqnarray*}
U(x,x+h]
&=& \int_{\bf R}Q(y,(x,x+h])\mu_0(dy)\\
&\le& \sup_{y\in(x,x+h]} Q(y,(x,x+h]).
\end{eqnarray*}
Therefore, it suffices to prove that
\begin{eqnarray}\label{renewal.boundedness}
Q(y,(x,x+h]) &\le& (A+h)/\varepsilon\delta
\end{eqnarray}
for all $y\in(x,x+h]$.
Given $X_0\in(x,x+h]$, consider the stopping time
$$
\tau=\min\{n\ge1:X_n>x+h\}.
$$
Since $X_\tau\wedge(x+h+A)-X_0 \le A+h$ with
probability $1$,
\begin{eqnarray*}
A+h &\ge& {\bf E}(X_\tau\wedge(x+h+A)-X_0)\\
&=& \sum_{n=1}^\infty {\bf E}
[X_n\wedge(x+h+A)-X_{n-1}\wedge(x+h+A)]{\bf I}\{\tau\ge n\}.
\end{eqnarray*}
Hence, the definition of $\tau$ implies
\begin{eqnarray*}
A+h &\ge& \sum_{n=1}^\infty
{\bf E}\{X_n\wedge(x+h+A)-X_{n-1}\wedge(x+h+A);\tau\ge n\}\\
&=& \sum_{n=1}^\infty
{\bf E}\{X_n\wedge(x+h+A)-X_{n-1}|\tau\ge n\}
{\bf P}\{\tau\ge n\}.
\end{eqnarray*}
The Markov property and condition
(\ref{sft.1}) yield
\begin{eqnarray*}
{\bf E}\{X_n\wedge(x+h+A)-X_{n-1}|\tau\ge n\}
&\ge& {\bf E}(\xi(X_{n-1})\wedge A)
\ge \varepsilon
\end{eqnarray*}
for all $n$. Therefore,
\begin{eqnarray*}
A+h &\ge& \varepsilon
\sum_{n=1}^\infty {\bf P}\{\tau\ge n\}
= \varepsilon {\bf E}\tau.
\end{eqnarray*}
So, the expected number of visits to the interval
$(x,x+h]$ till the first exit from
$(-\infty,x+h]$ does not exceed
$(A+h)/\varepsilon$, independently of
the initial state $X_0\in(x,x+h]$.
By condition (\ref{sft.2}),
after exit from $(-\infty,x+h]$
the chain is above the level $X_\tau$ forever
with probability at least $\delta$;
in particular, it does not visit the interval
$(x,x+h]$ any more.
With probability at most $1-\delta$
the chain visits this interval again, and so on.
Concluding, we get that the expected number of
visits to the interval $(x,x+h]$ cannot exceed the value
\begin{eqnarray*}
\frac{A+h}{\varepsilon} \sum_{n=0}^\infty (1-\delta)^n
&=& \frac{A+h}{\varepsilon\delta},
\end{eqnarray*}
and (\ref{renewal.boundedness}) is proved.
The proof of Theorem \ref{suff.for.tran} is complete.

The latter theorem yields the following

\begin{Corollary}
Let the family of jumps $\{\xi(x),x\in{\bf R}\}$
possess an integrable minorant with a positive mean,
that is, there exist a random variable $\zeta$
such that ${\bf E}\zeta>0$ and $\xi(x)\ge_{\rm st}\zeta$
for any $x\in{\bf R}$. Then
\begin{eqnarray*}
U(x,x+h] &\le& (A+h)A/\varepsilon^2
\end{eqnarray*}
for any $A>0$ such that
$\varepsilon\equiv{\bf E}(\zeta\wedge A)>0$;
in particular, {\rm(\ref{finite.bound.U})} holds.
\end{Corollary}

\proof. Consider the partial sums
$Z_n=\zeta_1+\ldots+\zeta_n$
of an independent copies of $\zeta$.
Denote the first ascending ladder epoch
by $\chi=\min\{n\ge1:Z_n>0\}$.
It is well known (see, for example, Theorem 2.3(c) in
[\ref{Aapq}, Ch. VIII]) that
\begin{eqnarray*}
{\bf P}\{Z_n>0\mbox{ for all }n\ge1\}
&=& 1/{\bf E}\chi.
\end{eqnarray*}
Since
\begin{eqnarray*}
{\bf P}\{X_n>x\mbox{ for all } n\ge 1|X_0=x\}
&\ge& {\bf P}\{Z_n>0\mbox{ for all }n\ge1\}
\end{eqnarray*}
by the minorisation condition,
the $\delta$ in Theorem \ref{suff.for.tran}
is at least $1/{\bf E}\chi$.
Taking into account the inequality
${\bf E}\chi\le A/\varepsilon$,
we get $\delta\ge \varepsilon/A$,
which implies the corollary conclusion.

If the chain $X$ has a non-negative jumps $\xi(x)\ge0$,
then the minorisation condition is equivalent to
the existence of a positive $A$ such that
\begin{eqnarray}\label{eps.non.neg}
\gamma \equiv \inf_{x\in{\bf R}}
{\bf P}\{\xi(x)>A\} &>& 0.
\end{eqnarray}
In that case one can choose $\zeta$ taking values
0 and $A$ with probabilities $1-\gamma$
and $\gamma$ respectively; then $\varepsilon\ge\gamma A$
and $U(x,x+h]\le(A+h)/\gamma^2A$.

%

We conclude with a counterexample demonstrating
that condition (\ref{majoriz})
in Theorem \ref{renewal.theorem} is essential;
the existence of integrable majorant cannot be
relaxed to the condition of the uniform
integrability of jumps.
We consider an integer valued chain $X$.
For any state $k\in{\bf Z}^+$, define the transition
probabilities in the following way: given $2^n-1\le k\le 2^{n+1}-2$,
\begin{eqnarray*}
p_{k,k} &=& 1/2,\\
p_{k,k+1} &=& 1/2-p_{k,2^{n+1}},\\
p_{k,2^{n+1}} &=& \frac{1}{(2^{n+1}-k)\ln(n+e^2)} \le \frac12.
\end{eqnarray*}
The corresponding jumps $\xi(k)$ converge
weakly as $k\to\infty$ to the Bernoulli distribution;
they are uniformly integrable.
But we can observe a concentration of relatively
large masses at points $2^{n+1}$;
the renewal measure at point $2^{n+1}$
is not less than, up to a positive constant,
\begin{eqnarray*}
\sum_{k=2^n-1}^{2^{n+1}-2} \frac{1}{(2^{n+1}-k)\ln(n+e^2)}
= \frac{1}{\ln(n+e^2)}\sum_{k=2}^{2^n+1} \frac{1}{k}
&\sim& \frac{n\ln 2}{\ln n}.
\end{eqnarray*}
Hence, there is no convergence $U\{k\}\to 1/{\bf E}\xi=2$
and the key renewal theorem
does not hold for the chain constructed.

This paper was mostly written while the author
was visiting the Boole Centre for Research in Informatics,
University College Cork,
thanks to the hospitality of Neil O'Connell
and financial support of Science Foundation
Ireland, grant no. SFI 04/RP1/I512.
Also this work was partially supported by
Russian Science Support Foundation.

\section*{\small REFERENCES}

\newcounter{bibcoun}
\begin{list}{\arabic{bibcoun}.}{\usecounter{bibcoun}\itemsep=0pt}
\small

\item\label{Alsmeyer}
Alsmeyer, G. (1994)
On the Markov renewal theorem.
{\it Stoch. Process. Appl.} {\bf 50}, 37--56.

\item\label{ADN}
Athreya, K. B., McDonald, D., and Ney, P. (1978).
Limit theorems for semi-Markov processes
and renewal theorem for Markov chains.
{\it Ann. Probab.} {\bf 6}, 788--797.

\item\label{Aapq}
Asmussen, S. (2003).
{\it Applied Probability and Queues},
Springer, New York.

\item\label{Borovkov}
Borovkov, A. A. (1999).
Asymptotically optimal solutions in a change-point problem.
{\it Theory Probab. Appl.} {\bf 43}, 539--561.

\item\label{BF}
Borovkov, A. A., and Foss, S. G. (1999).
Estimates for overshooting an arbitrary boundary
by a random walk and their applications.
{\it Theory Probab. Appl.} {\bf 44}, 231--253.

\item\label{CD}
Choquet, G., and Deny, J. (1960).
Sur l'\'equation de convolution $\mu=\mu*\sigma$.
{\it C. R. Acad. Sci. Paris S\'erie A} {\bf 250}, 799--801.

\item\label{F}
Feller, W. (1971).
{\it An Introduction to Probability Theory
and Its Applications},
Vol II, 2nd ed., John Wiley, New York.

\item\label{FO}
Feller, W., and Orey, S. (1961).
A renewal theorem.
{\it J. Math. Mech.} {\bf 10}, 619--624.

\item\label{Fuh}
Fuh, C.-D. (2004).
Uniform markov renewal theory and ruin
probabilities in Markov random walks.
{\it Ann. Appl. Probab.} {\bf 14}, 1202--1241.

\item\label{Heyde}
Heyde, C. C. (1967).
Asymptotic renewal results for a natural generalization
of classical renewal theory.
{\it J. Roy. Statist. Soc. Ser. B} {\bf 29}, 141--150.

\item\label{Horvath}
Horv\'ath, L. (1985).
A strong nonlinear renewal theorem with applications to
sequential analysis.
{\it Scand. J. Statist.} {\bf 12}, 271--280.

\item\label{Kesten}
Kesten, H. (1974).
Renewal theory for Markov chains.
{\it Ann. Probab.} {\bf 3}, 355--387.

\item\label{KW}
Kim, D.Y., and Woodroofe, M. (2006).
A non-linear Renewal Theorem with stationary
and slowly changing perturbations.
{\it I.M.S. Lecture Notes and Monograph Series},
{\bf 50}, 164--175.

\item\label{KP}
Kl\"uppelberg, C. and Pergamenchtchikov, S. (2003).
Renewal theory for functionals of a Markov chain
with compact state space.
{\it Ann. Probab.} {\bf 31}, 2270--2300.

\item\label{Korshunov}
Korshunov, D. A. (2001).
Limit theorems for general Markov chains.
{\it Sib. Math. J.} {\bf 42}, 301--316.

\item\label{LS1}
Lai, T. L., and Siegmund, D. O. (1977).
A non-linear renewal theory with applications
to sequential analysis, I.
{\it Ann. Statist.} {\bf 5}, 946-–954.

\item\label{LS2}
Lai, T. L., and Siegmund, D. O. (1979).
A non-linear renewal theory with applications
to sequential analysis, II.
{\it Ann. Statist.} {\bf 7}, 60-–76.

\item\label{M1975}
Maejima, M. (1975).
On local limit theorems and Blackwell's renewal theorem
for independent random variables.
{\it Ann. Inst. Statist. Math.}, {\bf 27}, 507-520.

\item\label{Tw}
Meyn, S. P., and Tweedie, R. L. (1993).
{\it Markov Chains and Stochastic Stability}.
Springer--Verlag, London.

\item\label{Nummelin}
Nummelin, E. (1978).
Uniform and ratio limit theorems for Markov renewal
and semi-regenerative processes on a general state space.
{\it Ann. I. H. Poincar\'e} {\bf 14}, 119--143.

\item\label{Raugi}
Raugi, A. (2004).
A general Choquet--Deny theorem for nilpotent groups.
{\it Ann. I. H. Poincar\'e} {\bf 40}, 677--683.

\item\label{Revuz}
Revuz, P. (1984).
{\it Markov Chains}.
North-Holland, Amsterdam.

\item\label{WW}
Wang, M., and Woodroofe, M. (1996).
A uniform renewal theorem.
{\it Sequential Anal.} {\bf 15}, 21--36.

\item\label{W1965}
Williamson, J. (1965).
Some renewal theorems for non-negative
independent random variables.
{\it Trans. Amer. Math. Soc.} {\bf 114}, 417-–445.

\item\label{Woodroofe}
Woodroofe, M. (1982).
{\it Nonlinear Renewal Theory in Sequential Analysis}.
CBMS-NSF Regional Conference Series in Applied
Mathematics, 39.
Society for Industrial and Applied Mathematics (SIAM),
Philadelphia, Pa.

\item\label{Zhang}
Zhang, C. H. (1988).
A non-linear renewal theory.
{\it Ann. Prob.} {\bf 16}, 793-–824.

\end{list}
\end{document}